\documentclass{article}
\usepackage{amsmath, epsfig, multicol}
\usepackage{amssymb}
\newtheorem{lemma}{Lemma}
\newtheorem{corollary}{Corollary}
\newtheorem{theorem}{Theorem}
\newtheorem{conjecture}{Conjecture}

\title{The Randi\'c index and the diameter of graphs}
\author{Yiting Yang
\thanks{Department of Mathematics, Zhejiang University, Hangzhou 310027, Zhejiang, PR China
({\tt yangyt@zju.edu.cn}).
This author was supported in part by  NSF grant
DMS 0701111.
}
  \and Linyuan Lu
\thanks{Department of Mathematics, University of South Carolina, Columbia, SC 29208,
({\tt lu@math.sc.edu}).
This author was supported in part by NSF grant
DMS 0701111 and DMS 1000475.
}}

\begin{document}
\maketitle
\begin{abstract}
  The {\it Randi\'c index} $R(G)$ of a graph $G$ is defined as the sum of
  $\frac{1}{\sqrt{d_ud_v}}$ over all edges $uv$ of $G$, where $d_u$
  and $d_v$ are the degrees of vertices $u$ and $v,$ respectively. Let $D(G)$ be the
  diameter of $G$ when $G$ is connected. Aouchiche-Hansen-Zheng \cite{Aouchiche}
  conjectured that among all connected graphs $G$ on $n$ vertices the
  path $P_n$ achieves the minimum values for both $R(G)/D(G)$ and
  $R(G)- D(G)$. We prove this conjecture completely.
  In fact, we prove a stronger theorem:  If $G$ is a connected graph, then
 $R(G)-\frac{1}{2}D(G)\geq \sqrt{2}-1$, with equality if and only if $G$ is a path with at
  least three vertices.
\end{abstract}

\section{Introduction}

In 1975, the chemist Milan Randi\'c \cite{Randic} proposed a
topological index $R$ under the name "branching index", suitable
for measuring the extent of branching of the carbon-atom skeleton
of saturated hydrocarbons. The branching index was renamed the
molecular connectivity index and is often referred to as the Randi\'c
index.

There is a good correlation between the Randi\'c index and several
physico-chemical properties of alkanes: boiling points, enthalpies
of formation,  chromatographic retention times, etc \cite{Hall, Kier1, Kier2}.

The {\it Randi\'c index} $R(G)$ of a graph $G=(V,E)$ is defined as
follows:
$$R(G)=\sum_{uv\in E} \frac{1}{\sqrt{d_u d_v}}.$$
Here $d_u$ and $d_v$ are the degrees of vertices $u$ and $v$, respectively.

From a mathematical point of view, the first question to be asked
is what are the minimum and maximum values of Randi\'c index in  various
classes of graphs, and which graphs in these classes of
graphs have an extremal (minimum or maximum) Randi\'c index.
Erd\H{o}s and Bollob\'as \cite{Bollobas} first considered such
problems. They proved that the star minimizes the  Randi\' c
index among all the graphs without isolated vertices on fixed
number of vertices. After that a lot of extremal results on the
Randi\'c index were published.

It turns out that the Randi\'c index is also related to some typical
graph parameters such as:  diameter, radius, average distance,
girth, chromatic number, and eigenvalues of the adjacent matrices
 \cite{Aroujo1,Aroujo2,Li1}. Some conjectures on them are still
open \cite{Fajtlowicz, Aouchiche, Li2}.

Aouchiche-Hansen-Zheng \cite{Aouchiche} posed the following conjecture
on the diameter and the Randi\'c index.

\begin{conjecture}\label{cj}
If $G$ is a connected graph  of order $n\geq 3$, 
then  the Randi\'c index $R(G)$ and the diameter $D(G)$ satisfy
$$R(G)-D(G)\geq \sqrt{2}-\frac{n+1}{2}$$ and
$$\frac{R(G)}{D(G)}\geq \frac{n-3+2\sqrt{2}}{2n-2},$$
with equalities if and only if $G \cong P_n.$
\end{conjecture}

Li and Shi \cite{lishi} proved this conjecture in  some special cases.
Namely, if $G$ is  a connected graph of order $n$ with minimum
degree at least $5$, then $$R(G)-D(G)\geq  \sqrt{2}-\frac{n+1}{2}.$$
If $\delta(G)\geq \frac{n}{5}$, then $$\frac{R(G)}{D(G)}\geq \frac{n-3+2\sqrt{2}}{2n-2}.$$

In this paper we settle the conjecture completely.
In fact, we prove the following stronger theorem.

\begin{theorem}\label{main}
 If $G$ is a connected graph with at least three vertices,
then we have
$$R(G)-\frac{1}{2}D(G)\geq \sqrt{2}-1.$$
Equality holds if and only if $G\cong P_n$ for $n\geq 3$.
\end{theorem}

\begin{corollary}
If $G$ is a connected graph  of order $n\geq 3$, 
then  the Randi\'c index $R(G)$ and the diameter $D(G)$ satisfy
$$R(G)-D(G)\geq \sqrt{2}-\frac{n+1}{2}$$ and
$$\frac{R(G)}{D(G)}\geq \frac{n-3+2\sqrt{2}}{2n-2},$$
with equalities if and only if $G \cong P_n.$
\end{corollary}
\noindent
{\bf Proof.} 
Noticing that $D(G)\leq n-1$, we have

$$R(G)-D(G)=R(G)-\frac{D(G)}{2}-\frac{D(G)}{2} \geq
\sqrt{2}-1-\frac{n-1}{2}=\sqrt{2}-\frac{n+1}{2}$$


and

$$R(G)-\frac{D(G)}{2}\geq \sqrt{2}-1 \Rightarrow
\frac{R(G)}{D(G)}\geq \frac{1}{2}+\frac{\sqrt{2}-1}{D(G)}\geq
\frac{n-3+2\sqrt{2}}{2n-2}.$$

\hfill $\square$

The paper is organized as follows. In section 2, we prove several
useful lemmas. Our main idea is to capture the change of
the Randi\'c index when we simplify a graph.
The proof of the main theorem is presented in section 3.

\section{Lemmas on vertex deletion and edge deletion}

For any vertex $v$, let $\Gamma(v)$ denote the set of all neighbors of $v$ and
$\Gamma^*(v)$ denote the set of all  non-leaf neighbors of $v,$
i. e.,
$$\Gamma(v)=\{u \colon uv\in E(G)\}\ \  \mbox{and}\ \  \Gamma^*(v)=\{u \colon uv\in E(G)\ \mbox{ and }\  d_u\geq 2\}.$$
We also let $N(v)=\Gamma(v)\cup\{v\}$ and $N^*(v)= \Gamma^*(v)\cup \{v\}$.
Throughout $d_u$ will be degree with respect to $G$, unless other
graphs are considered. 

We have the following Lemma.
\begin{lemma}\label{l0}
  If $G$ is a connected graph on the vertex set $ \{1,2,\ldots, n\}$, then we have
$$R(G)\geq \frac{\sum_{i=1}^n \sqrt{d_i}}{2 \sqrt{\Delta}}.$$
Here $d_1,\dots, d_n$ are degrees of $G$ and $\Delta$ is the maximum degree.
\end{lemma}
{\bf Proof:} We have
\begin{eqnarray*}
R(G)&=&\sum_{ij \in E(G)}\frac{1}{\sqrt{d_id_j}}\\
&=&\frac{1}{2}\sum_{i=1}^n\sum_{j \in \Gamma(i)}\frac{1}{\sqrt{d_id_j}}\\
&\geq& \frac{1}{2}\sum_{i=1}^n\sum_{j \in \Gamma(i)}\frac{1}{\sqrt{d_i\Delta}}\\
&\geq& \frac{1}{2}\sum_{i=1}^n\frac{d_i}{\sqrt{d_i\Delta}}\\
&=&\frac{\sum_{i=1}^n \sqrt{d_i}}{2 \sqrt{\Delta}}.
\end{eqnarray*}
The proof of this lemma is finished. \hfill $\square$

Let $G - v$ be the induced subgraph obtained by
deleting the vertex $v$ from $G$. Let $G - uv$ be the
spanning subgraph obtained by deleting the edge $uv$ from $G$.

If $G$ is connected, then $D(G)$ is the diameter of $G$ as defined early.
We extend the function $D(G)$ to disconnected graphs as follows.
If $G$ is disconnected, then $D(G)$ is defined to be the
maximum among diameters of all the connected components of $G$.
A vertex $v$ is said to be {\it essential} (to $D(G)$) if $D(G - v)<
D(G)$; it is not essential otherwise. Thus a vertex $v$ is
essential if and only if every shortest path between any two vertices at
distance $D(G)$ passes through $v$.

An edge is {\it essential} if its two
endpoints are essential. A path is {\it essential} if all edges of
this path are essential.

In general, $\Gamma(v)$ is not an independent set. Let
$G|_{\Gamma(v)}$ be the induced subgraph of $G$ on $\Gamma(v).$ We
have the following lemma.

\begin{lemma}\label{l3}
Given an orientation of the edges of $G|_{\Gamma(v)}$, for any two
vertices $u$ and $x$ in $G$,
we define
\[\epsilon^u_x=\left\{
  \begin{array}{cc}
    1, & \mbox{ if $\overrightarrow{ux}$ is a directed edge of $G|_{\Gamma(v)}$;}\\
    0, & \mbox{ otherwise.}
  \end{array}
\right.
\]
If for any $u\in\Gamma^*(v),$
\begin{equation}
  \label{eq:1}
 \frac{1}{d_u-1}
\sum_{x\in \Gamma(u)\setminus
\{v\}}\frac{1}{\sqrt{d_x-\epsilon^u_x}} \leq \frac{2}{\sqrt{d_v}},
\end{equation}
then we have
$$R(G)> R(G - v).$$
\end{lemma}
{\bf Proof:} When the vertex $v$ is deleted, all edges incident
to $v$ are also deleted.  For any vertex $u$, if $u\in \Gamma(v)$, the
degree of $u$ decreases by one; if $u\not\in N(v)$,
the degree of $u$ remains the same.

Let us consider $R(G)-R(G - v)$. For most edges $xy$ in $G$, the
contribution of $\frac{1}{\sqrt{d_xd_y}}$ to $R(G)-R(G - v)$ is
canceled out unless one of $x$ and $y$ is in $N(v)$. There are
three types of edges.

\noindent
{\bf Type I:}  $x=v$ and $y=u\in \Gamma(v)$. The contribution of this type of edge
 to $R(G)-R(G - v)$ is
$$\sum_{u\in\Gamma(v)}\frac{1}{\sqrt{d_vd_u}}\geq \sum_{u\in\Gamma^*(v)}\frac{1}{\sqrt{d_vd_u}} .$$

\noindent {\bf Type II:}  $y=u\in \Gamma^*(v)$ and $x\in \Gamma(u)\setminus N(v)$.
The contribution of this type of edge
 to $R(G)-R(G - v)$ is
 \begin{eqnarray*}
&&\hspace*{-1cm}\sum_{u\in\Gamma^*(v)}\left(\frac{1}{\sqrt{d_u}}-\frac{1}{\sqrt{d_u-1}}\right)
\sum_{x\in \Gamma(u)\setminus N(v)}
\frac{1}{\sqrt{d_x}} \\
&=& \sum_{u\in\Gamma^*(v)}\left(\frac{1}{\sqrt{d_u}}-\frac{1}{\sqrt{d_u-1}}\right)
\sum_{x\in \Gamma(u)\setminus N(v)}
\frac{1}{\sqrt{d_x-\epsilon_x^u}}
    \end{eqnarray*}
since $\epsilon_x^u=0$ in this case.

\noindent
{\bf Type III:}  $y=u\in \Gamma^*(v)$, $x\in  \Gamma^*(v)$, and
$\overrightarrow{ux}$ is a directed edge of $G|_{\Gamma(v)}$.
Note that
\begin{eqnarray*}
&&\hspace*{-1cm}
\frac{1}{\sqrt{d_ud_x}}-\frac{1}{\sqrt{(d_u-1)(d_x-1)}} \\
&=&
\frac{1}{\sqrt{d_u}}\left(\frac{1}{\sqrt{d_x}}-\frac{1}{\sqrt{d_x-1}}\right)
+
\frac{1}{\sqrt{d_x-1}}\left(\frac{1}{\sqrt{d_u}}-\frac{1}{\sqrt{d_u-1}}\right)\\
&=&
\frac{1}{\sqrt{d_u-\epsilon_u^x}}\left(\frac{1}{\sqrt{d_x}}-\frac{1}{\sqrt{d_x-1}}\right)
+
\frac{1}{\sqrt{d_x-\epsilon_x^u}}\left(\frac{1}{\sqrt{d_u}}-\frac{1}{\sqrt{d_u-1}}\right),
\end{eqnarray*}
since $\epsilon_u^x=0$ and $\epsilon_x^u=1$. The above expression is symmetric
with respect to $u$ and $x$. Thus, the contribution of this type of edge
 to $R(G)-R(G - v)$ is
 \begin{eqnarray*}
&&\hspace*{-1cm}\frac{1}{2} \sum_{u\in \Gamma^*(v), x\in \Gamma(u)\cap \Gamma(v)}
\frac{1}{\sqrt{d_u-\epsilon_u^x}}\left(\frac{1}{\sqrt{d_x}}-\frac{1}{\sqrt{d_x-1}}\right)
+
\frac{1}{\sqrt{d_x-\epsilon_x^u}}\left(\frac{1}{\sqrt{d_u}}-\frac{1}{\sqrt{d_u-1}}\right)\\
&=&\sum_{u\in \Gamma^*(v), x\in \Gamma(u)\cap \Gamma(v)}\frac{1}{\sqrt{d_x-\epsilon_x^u}}\left(\frac{1}{\sqrt{d_u}}-\frac{1}{\sqrt{d_u-1}}\right)\\
&=&\sum_{u\in \Gamma^*(v)}\left(\frac{1}{\sqrt{d_u}}-\frac{1}{\sqrt{d_u-1}}\right)
\sum_{ x\in \Gamma(u)\cap \Gamma(v)} \frac{1}{\sqrt{d_x-\epsilon_x^u}}.
 \end{eqnarray*}

Summing up the contribution of three types of edges,
we have
\begin{eqnarray*}
  R(G)-R(G - v)&\geq&\sum_{u\in\Gamma^*(v)}\frac{1}{\sqrt{d_vd_u}} + \sum_{u\in\Gamma^*(v)}\left(\frac{1}{\sqrt{d_u}}-\frac{1}{\sqrt{d_u-1}}\right)
\sum_{x\in \Gamma(u)\setminus N(v)}
\frac{1}{\sqrt{d_x-\epsilon_x^u}}\\
&&+\sum_{u\in \Gamma^*(v)}\left(\frac{1}{\sqrt{d_u}}-\frac{1}{\sqrt{d_u-1}}\right)
\sum_{ x\in \Gamma(u)\cap \Gamma(v)} \frac{1}{\sqrt{d_x-\epsilon_x^u}}.\\
&=& \sum_{u\in\Gamma^*(v)}\left[\frac{1}{\sqrt{d_vd_u}}
-\left( \frac{1}{\sqrt{d_u-1}}
-\frac{1}{\sqrt{d_u}}\right)\sum_{x\in\Gamma(u)\setminus\{v\}}\frac{1}{\sqrt{d_x-\epsilon_x^u}}\right].
\end{eqnarray*}
Now we apply the assumption (\ref{eq:1}).
\begin{eqnarray*}
  R(G)-R(G - v)
&\geq& \sum_{u\in \Gamma^{*}(v)} \left[ \frac{1}{\sqrt{d_vd_u}}-
\left( \frac{1}{\sqrt{d_u-1}}
-\frac{1}{\sqrt{d_u}}\right) \frac{2(d_u-1)}{\sqrt{d_v}}
\right] \\
&=&\sum_{u\in \Gamma^{*}(v)}\frac{1}{\sqrt{d_vd_u}} \left(1-
\frac{2\sqrt{d_u-1}}{\sqrt{d_u}+\sqrt{d_u-1}}
\right) \\
&=&\sum_{u\in \Gamma^{*}(v)}\frac{(\sqrt{d_u}-\sqrt{d_u-1})^2}{\sqrt{d_vd_u}}\\
&>&0. \hspace*{3in} \square
\end{eqnarray*}

Inequality (\ref{eq:1}) is called the {\it deletion condition} for the
vertex $v$. To check the deletion condition, we need to specify
an orientation of  the edges of $G\mid_{\Gamma(v)}$. We can
relax this condition as follows.

Let $d^*_x=d_x-1$ if $d_x\geq 2$ and $d^*_x=d_x$ if $d_x=1$.
Note for any orientation of the edges of $G\mid_{\Gamma(v)}$
$$d_x-\epsilon_x^u\geq d^*_x.$$
We have the following corollary.
\begin{corollary} \label{cor2}
If for any $u\in \Gamma^*(v)$,
\begin{equation}
  \label{eq:2}
 \frac{1}{d_u-1}
\sum_{x\in \Gamma(u)\setminus \{v\}}\frac{1}{\sqrt{d_x^*}} \leq
\frac{2}{\sqrt{d_v}},
\end{equation}
then we have $$R(G)> R(G - v).$$
\end{corollary}

Inequality (\ref{eq:2}) is called the {\it weak deletion condition} for the
vertex $v$.

\begin{corollary} \label{cor3}
  If $d_v\leq 4$,  then we have
$$R(G)>R(G - v).$$
\end{corollary}
{\bf Proof:} It suffices to show that $v$
satisfies the weak deletion condition.  If $\Gamma^*(v)=\emptyset$, then the
weak deletion condition is satisfied automatically.
If $u\in \Gamma^*(v)$ and $x\in \Gamma(u)\setminus\{v\}$,  then we have
$$d_x^*\geq 1.$$
Thus,
 \[\frac{1}{d_u-1}
\sum_{x\in \Gamma^*(u)\setminus
\{v\}}\frac{1}{\sqrt{d_x^*}}\leq
\frac{1}{d_u-1}\sum_{x\in \Gamma^*(u)\setminus \{v\}}1\leq 1
\leq\frac{2}{\sqrt{d_v}}.\]
Applying Corollary \ref{cor2}, we get
$$R(G)> R(G - v).$$
\hfill $\square$

\begin{lemma}\label{l4}
If $G$ is a connected graph,  then there exists an induced connected subgraph $G'$
satisfying the following conditions.
  \begin{enumerate}
  \item $R(G)\geq R(G')$.
  \item $D(G)\leq D(G')$.
  \item Every non-essential vertex in $G'$ has degree at least $9$.
  \item $R(G')=R(G)$ holds if and only if $G'=G$ and every non-essential vertex
in $G$ has degree at least $9$.
  \end{enumerate}

\end{lemma}
{\bf Proof:} Suppose that $G$ contains a vertex $v$ with $d_v\leq 4$.
If $v$ is not essential, then we can remove $v$ from $G$ and consider
$G - v$ instead (by Corollary \ref{cor3}). Repeatedly find a non-essential vertex $v$ with degree at most $4$
and delete it until no such $v$ is found.

From now on, we assume every non-essential vertex has degree at
least $5$. Let $v$ be a non-essential vertex with minimum degree
$\delta\leq 8$. We claim $$R(G)> R(G - v).$$ There are five cases.

\noindent
{\bf Case I:}  The vertex $v$ has one neighbor $u_1$ with degree $1$, and $u_1$ is
essential. Any path containing $u_1$ contains $v$. This contradicts with the assumption that $v$ is not essential.

\noindent
{\bf Case II:} The vertex
 $v$ has two neighbors $u_1$ and $u_2$ with degrees  $2$, and
both $u_1$ and $u_2$ are essential vertices. Since $v$ is not
essential, there exists a shortest path $P$ (of length $D(G)$)
which  does not contain $v$.
The path $P$ passes through $u_1$ and $u_2$. The degrees of $u_1$ and $u_2$
in $P$ are at most $1$. So $u_1$ and $u_2$ must be the two endpoints
of $P$. In this case, we must have $D(G)=d(u_1,u_2)\leq 2$.

If $D(G)=1$, then $G$ is a complete graph. We have
$R(G)=R(G - v)+\frac{1}{2}>R(G - v)$.

Now assume $D(G)=2$. Since $d_v=\delta\geq 5$,
 $\Gamma(v)$ contains  a vertex $u$ which is not on the path $P$.
We have $d(u, u_i)=2$ for $i=1,2$. We can delete $u_1$ or $u_2$ without
decreasing $D(G)$. Contradiction!

\noindent
{\bf Case III:}  Every neighbor of $v$ has degree at least $3$, and no leaf
lies within the distance $2$ from $v$. For any $u\in \Gamma(v)$
with degree at least $3$ and $x\in \Gamma(u)\setminus \{v\}$, we
have
$$d_x^*\geq 2.$$
We have
 \[\frac{1}{d_u-1}
\sum_{x\in \Gamma(u)\setminus \{v\}}\frac{1}{\sqrt{d_x^*}} \leq
\frac{1}{\sqrt{2}} \leq\frac{2}{\sqrt{d_v}},\] which holds for
$d_v\leq 8$. The weak deletion condition (\ref{eq:2}) is satisfied. By
Corollary \ref{cor2}, we have
$$R(G)> R(G - v).$$

\noindent
{\bf Case IV:}
All neighbors of $v$ except $u_1$  have degree at least $3$ while $u_1$ has degree $2$;
 no leaf lies within the distance $2$ from $v$. In this case, we verify
 the deletion condition (\ref{eq:1}).
Orient the edges of $G|_{\Gamma^*(v)}$ so that the edge incidents
to $u_1$ leave $u_1$. For any $u\in \Gamma(v)$ and $x\in
\Gamma^*(u)\setminus \{v\}$, it is clear that
$$d_x-\epsilon^u_x\geq 2.$$ Similarly, the condition (\ref{eq:1})
is satisfied. By Lemma \ref{l3}, we have
$$R(G)> R(G - v).$$

\noindent
{\bf Case V:}
 There is a leaf $x$ with $d(v,x)=2$,  and $x$ is essential. Let
$u$ be the only neighbor of $x$. Clearly, $u\in \Gamma^*(v)$. Since $x$ is
essential, then $u$ must be essential as well. We verify the weak deletion
condition (\ref{eq:2}) for $u$.

If $d_u=2$, then $v$ is also essential. Contradiction!
Suppose that $u$ has a neighbor $w$ with $d_w<\delta$ ( $w\not=x$).
The vertex $w$ must be essential. Since $d_v=\delta>d_w$, there
is a vertex $y\in \Gamma(v)\setminus \Gamma(w)$. Suppose that $P$ is a shortest path of length $D(G)$
containing $x, u, w$. Replace the segment $x-u-w$ by the shortest path from
$y$ to $w$. Call this path $P'$. The path $P'$ is also a shortest path with
length at least $D(G)$, and $P'$ dose not  contain  $x$. This contradicts with the assumption that $x$ is essential.

Suppose $d_u\geq 3$, and every neighbor $w$ of $u$ other than $x$
satisfies $d_w\geq \delta$.
We have
\begin{eqnarray*}
&& \hspace*{-2cm} \frac{1}{d_u-1} \sum_{x\in \Gamma^*(u)\setminus
\{v\}}\frac{1}{\sqrt{d_x^*}}\\
&\leq&  \frac{1}{d_u-1} \left(1 + (d_u-2)\frac{1}{\sqrt{\delta-1}}\right)\\
&=& \frac{1}{\sqrt{\delta-1}} + \frac{1}{d_u-1}\left(1-\frac{1}{\sqrt{\delta-1}}\right)\\
&\leq& \frac{1}{\sqrt{\delta-1}} + \frac{1}{2}\left(1-\frac{1}{\sqrt{\delta-1}}\right)\\
&\leq& \frac{1}{2}\left(1+\frac{1}{\sqrt{\delta-1}}\right)\\
&<& \frac{2}{\sqrt{\delta}}.
\end{eqnarray*}
The last inequality holds for $\delta\leq 8$. Thus $R(G)\geq R(G - v)$.

For all five cases, we can delete a non-essential vertex $v$ with $d_v\leq 8$
 such that
$$R(G)> R(G - v).$$
Repeat this process until every non-essential vertex has degree
at least $9$.
\hfill $\square$

A vertex $v$ is  a {\em local-minimum-vertex} if the following two conditions
are satisfied.
\begin{enumerate}
\item The vertex $v$ is not essential for $G$.
\item If $u$ is a  non-essential vertex with $d(u,v)\leq 2$, then $d_u\geq d_v$.
\end{enumerate}

\begin{lemma}\label{local}
Suppose $v$ is a local-minimum-vertex with degree $d_v\geq 3$. If
$R(G)\leq R(G - v)$, then there exist two vertices $w$ and $y$
satisfying
\begin{enumerate}
\item  $vw$ and $wy$ are edges of $G$.
\item  $d_w<d_v$ and $d_y<d_v$. Consequently, $wy$ is
an essential edge of $G$.
\end{enumerate}
\end{lemma}
{\bf Proof:} For any $u\in \Gamma(v)$ with $d_u\geq d_v$,
we claim that $\Gamma(u)$ can contain at most two essential vertices.

Otherwise, say that $\Gamma(u)$ contains three essential
vertices $x$, $y$, and $z$. Choose a shortest path $P$ connecting
two vertices of distance $D(G)$. By the definition of essential vertices,
all $x$, $y$, and $z$ are on the path $P$. Since $x$, $y$, $z\in \Gamma(u)$,
 $x$, $y$, and $z$ must be
adjacent on $P$. Without loss of generality, we assume $d(x,z)=2$.
We can replace $y$ by $u$ and obtain a new path $P'$ from $P$.
This contradicts with the assumption that $y$ is also essential.

Since $R(G)\leq R(G - v)$, the weak deletion condition (\ref{eq:2})
is violated for some $u$. There are three cases.

\noindent
{\bf Case I:} For any $x\in \Gamma(u)\setminus \{v\}$, $d_x\geq d_v$.
In this case, we have
$$d^*_x=d_x-1\geq d_v-1.$$
Thus,
\begin{eqnarray*}
 \frac{1}{d_u-1}
\sum_{x\in \Gamma(u)\setminus \{v\}}\frac{1}{\sqrt{d_x^*}}
&\leq&  \frac{1}{d_u-1}
\sum_{x\in \Gamma(u)\setminus \{v\}} \frac{1}{\sqrt{d_v-1}}\\
&\leq&  \frac{1}{\sqrt{d_v-1}}\\
&<& \frac{2}{\sqrt{d_v}},
\end{eqnarray*}
where the lase step holds for $d_v\geq 2$.  The weak deletion condition
(\ref{eq:2}) is satisfied. Contradiction!

\noindent
{\bf Case II:}  $d_u<d_v$. By Case I, we have a vertex
$x\in \Gamma(u)\setminus \{v\}$, $d_x< d_v$.
Choose $w=u$ and $y=x$. We are done.

\noindent
{\bf Case III:}
$d_u\geq d_v$. Note that $\Gamma(u)$ can contain at most $2$
  essential vertices. Let $y_1$ and $y_2$ be the possible two
  essential vertices. If $x\in \Gamma(u)\setminus \{v,y_1,y_2\}$,
   then by the definition of local-minimum-vertex, we have
$$d^*_x=d_x-1\geq d_v-1.$$
We  bound $d^*_{y_1}$ and $d^*_{y_1}$  by $2$.
We get
\begin{eqnarray*}
 \frac{1}{d_u-1}
\sum_{x\in \Gamma(u)\setminus \{v\}}\frac{1}{\sqrt{d_x^*}}
&\leq&  \frac{1}{d_u-1}\left[ 2 + \frac{d_u-3}{\sqrt{d_v-1}}\right]
\\
&\leq&  \frac{1}{\sqrt{d_v-1}} +\frac{2-\frac{2}{\sqrt{d_v-1}}}{d_u-1}  \\
&\leq&  \frac{1}{\sqrt{d_v-1}} +\frac{2-\frac{2}{\sqrt{d_v-1}}}{d_v-1}  \\
&<& \frac{2}{\sqrt{d_v}},
\end{eqnarray*}
where the lase step holds for $d_v\geq 3$. Contradiction!





Only Case II is possible. There are two essential vertices $y$ and
$w$ satisfying all the conditions.  \hfill $\square$


\begin{lemma}\label{cut}
If $uv$ is a non-leaf edge, then we have
$$R(G)> R(G - uv)-\frac{1}{2}.$$
\end{lemma}
{\bf Proof:} We have
\begin{eqnarray*}
  R(G)-R(G - uv) &=& \frac{1}{\sqrt{d_ud_v}} -
 \sum_{x\in \Gamma(u)\setminus \{v\}}\frac{1}{\sqrt{d_x}}
\left(\frac{1}{\sqrt{d_u-1}}-\frac{1}{\sqrt{d_u}}\right)\\
&\ & - \sum_{y\in \Gamma(v)\setminus \{u\}}\frac{1}{\sqrt{d_y}}
\left(\frac{1}{\sqrt{d_v-1}}-\frac{1}{\sqrt{d_v}}\right)\\
&\geq& \frac{1}{\sqrt{d_ud_v}} -
(d_u-1)\left(\frac{1}{\sqrt{d_u-1}}-\frac{1}{\sqrt{d_u}}\right)\\
&\ &-(d_v-1)\left(\frac{1}{\sqrt{d_v-1}}-\frac{1}{\sqrt{d_v}}\right)\\
&=& \frac{1}{\sqrt{d_ud_v}} - \frac{\sqrt{d_u-1}}{\sqrt{d_u}(\sqrt{d_u}+\sqrt{d_u-1})}
-\frac{\sqrt{d_v-1}}{\sqrt{d_v}(\sqrt{d_v}+\sqrt{d_v-1})}\\
&>& \frac{1}{\sqrt{d_ud_v}}-\frac{1}{2\sqrt{d_u}}-\frac{1}{2\sqrt{d_v}}\\
&=& \frac{1}{\sqrt{2d_ud_v}} +
\frac{1}{2}\left( 1-\frac{1}{\sqrt{d_u}}\right)\left( 1-\frac{1}{\sqrt{d_v}}\right) -\frac{1}{2}\\
&>&-\frac{1}{2}.
\end{eqnarray*}
\hfill $\square$

\begin{corollary}\label{cut1}
  Suppose that  $uv$ is not a cut edge of $G$. If both $u$ and $v$ are essential, then
$$R(G)-\frac{1}{2}D(G)>R(G - uv)-\frac{1}{2}D(G - uv).$$
\end{corollary}


\begin{lemma}\label{attach}
Let $u$ be a cut vertex of $G$. Suppose that $G$ has a decomposition $G=G_1\cup G_2$
satisfying $G_1\cap G_2=\{u\}$,  $|G_2|\geq 8$, and  $|\Gamma_u\cap V(G_1)|=2$
(see Figure \ref{fig:decomp}).  If $u$ reaches the minimum degree in $G_2$, then
we have
$$R(G)> R(G_1).$$
\end{lemma}

\begin{center}
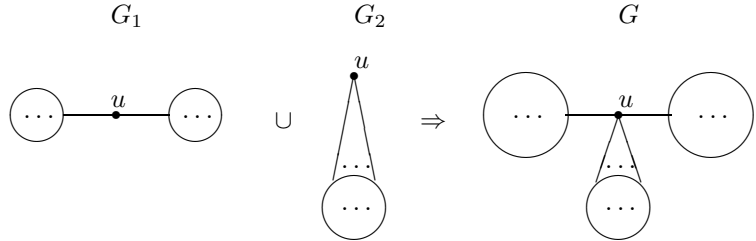
\begin{figure}[hbt]

\begin{picture}(290,100)
\put(20,45){\circle{20}}

\put(15,44){$\ldots$}

\put(30,45){\line(1,0){40}}

\put(50,45){\circle*{3}}

\put(48,48){$u$}

\put(80,45){\circle{20}}

\put(75,44){$\ldots$}

\put(48,80){$G_1$}

\put(110,40){$\cup$}

\put(140,80){$G_2$}

\put(140,10){\circle{25}}

\put(135,9){$\ldots$}

\put(140,60){\circle*{3}}

\put(140,63){$u$}

\put(140,60){\line(1,-5){7.9}}

\put(140,60){\line(-1,-5){7.9}}

\put(135,25){$\ldots$}

\put(165,40){$\Rightarrow$}

\put(205,45){\circle{30}}

\put(200,44){$\ldots$}

\put(220,45){\line(1,0){40}}

\put(240,80){$G$}

\put(240,48){$u$}

\put(240,45){\circle*{3}}

\put(240,45){\line(1,-3){8.5}}

\put(240,45){\line(-1,-3){8.5}}

\put(240,10){\circle{25}}

\put(235,25){$\ldots$}

\put(235,9){$\ldots$}

\put(275,45){\circle{30}}

\put(270,44){$\ldots$}

\end{picture}
\caption{$G=G_1\cup G_2$}
\label{fig:decomp}
\end{figure}
\end{center}

{\bf Proof:} Let $u_1$ and $u_2$ be the two adjacent vertices of $u$ in $G_1$
and $v_1,\ldots, v_k$ be the adjacent vertices of $u$ in $G_2.$
Let $N(v_i)$ be the set of neighbors of $v_i$ in $G_2$.  We
have
\begin{align*}
R(G)&\geq
R(G_1)+R(G_2)-\left(\frac{1}{\sqrt{2}}-\frac{1}{\sqrt{k+2}}\right)\left(\frac{1}{\sqrt{d_{u_1}}}+\frac{1}{\sqrt{d_{u_2}}}\right)\\
&\
-\sum_{i=1}^{k}\left(\frac{1}{\sqrt{d_{v_i}}}-\frac{1}{\sqrt{d_{v_i}+1}}\right)\sum_{x\in
N(v_i)}\frac{1}{\sqrt{d_x}}\\
&\geq
R(G_1)+R(G_2)-\frac{1}{\sqrt{2}}\left(\frac{1}{\sqrt{d_{u_1}}}+\frac{1}{\sqrt{d_{u_2}}}\right)\\
&\
-\sum_{i=1}^k\frac{1}{(\sqrt{d_{v_i}}+\sqrt{d_{v_i}+1})\sqrt{d_{v_i}}\sqrt{d_{v_i}+1}}\cdot \frac{d_{v_i}}{\sqrt{k}}\\
& >R(G_1)+R(G_2)-\frac{1}{\sqrt{2}}\left(\frac{1}{\sqrt{d_{u_1}}}+\frac{1}{\sqrt{d_{u_2}}}\right)
-\sum_{i=1}^{k}\frac{1}{2\sqrt{k}\sqrt{d_{v_i}}}\\
&\geq
R(G_1)+R(G_2)-\frac{1}{\sqrt{2}}\left(\frac{1}{\sqrt{1}}+\frac{1}{\sqrt{1}}\right)-\frac{k}{2\sqrt{k}\sqrt{k}}\\
&\geq R(G_1)+R(G_2)-\sqrt{2}-\frac{1}{2}>R(G_1),
\end{align*}
where the last inequality hold for $R(G_2)\geq \sqrt{|G_2|-1}\geq
\sqrt{7}.$

\hfill $\square$

\begin{lemma}\label{subdivision}
 For any edge $uv$ of $G$, let $G_{u\cdot v}$ be the graph obtained by
subdividing the edge $uv$, (i.e., by replacing the edge $uv$ by a path
of length $2$.) We have the following statements.
\begin{enumerate}
\item If $d_u=2$ or $d_v=2$, then $R(G_{u\cdot v})=R(G)+\frac{1}{2}$.
\item If $d_u>2$ and $d_v>2$, then $R(G_{u\cdot v})<R(G)+\frac{1}{2}$.
\item If $d_u=1$ and $d_v>2$, then  $R(G_{u\cdot v})>R(G)+\frac{1}{2}$.
\item If $d_u>2$ and $d_v=1$, then  $R(G_{u\cdot v})>R(G)+\frac{1}{2}$.
\end{enumerate}
\end{lemma}
{\bf Proof:}
We have
\begin{eqnarray*}
 R(G_{u\cdot v})-R(G)&=& \frac{1}{\sqrt{2d_u}}+ \frac{1}{\sqrt{2d_v}}
-\frac{1}{\sqrt{d_ud_v}} \\
&=& \frac{1}{2} - \left(\frac{1}{\sqrt{2}}-\frac{1}{\sqrt{d_u}}\right)
\left(\frac{1}{\sqrt{2}}-\frac{1}{\sqrt{d_v}}\right).
\end{eqnarray*}
It is easy to verify all cases. \hfill $\square$

\section{Proof of  main theorem}

{\bf Proof of Theorem \ref{main}:} For any graph $G$, we define
$f(G)=R(G)-\frac{D(G)}{2}$.  Note that $f(P_n)=\sqrt{2}-1$ for
$n\geq 3$. We need show that
\begin{equation}\label{minequlity}
f(G)> \sqrt{2}-1
\end{equation}
for any connected graph $G\not=P_n$ ($n\geq 3$).

Suppose that there is such a  graph $G$ ($\not=P_n$) satisfying
$$f(G)\leq\sqrt{2}-1.$$
Let $G$ be such a graph with the smallest number of vertices.
(If there are several such graphs with the same number of vertices,
pick the one with minimum number of edges.)
It is easy to check that $G$ is connected and has at least 3 vertices.

By Lemma \ref{l4}, every non-essential vertex of $G$ has degree at
least 9. By Corollary \ref{cut1}, every essential edge is an
edge-cut of $G$. By Lemma \ref{attach}, if there are two essential
edges $uv$ and $vw$, then $d_v=2$. Therefore $G$ is the graph
consists of several blocks which are linked by essential paths
(see Figure \ref{figure2}). A {\em block} $B$ is an induced
connected subgraph of $G$ which contains no essential edges of
$G$. By Lemma \ref{subdivision}, the length of each essential path
is either $1$ or $2$.

\begin{center}
  \begin{figure}[hbt]
    \centering
\begin{picture}(360, 100)

\put(10,40){\line(1,0){30}}

\put(10,60){\line(1,0){30}}

\put(10,40){\line(0,1){20}}

\put(40,40){\line(0,1){20}}

\put(40,50){\line(1,0){20}}

\put(80,50){\line(1,0){20}}

\put(100,30){\line(0,1){40}}

\put(260,30){\line(0,1){40}}

\put(100,30){\line(1,0){160}}

\put(179,15){$B$}

\put(100,70){\line(1,0){160}}

\put(260,50){\line(1,0){20}}

\put(300,50){\line(1,0){20}}

\put(320,40){\line(0,1){20}}

\put(320,40){\line(1,0){30}}

\put(350,40){\line(0,1){20}}

\put(320,60){\line(1,0){30}}

\put(100,50){\line(2,1){20}}

\put(100,50){\line(2,-1){20}}

\put(140,50){\line(-2,1){20}}

\put(140,50){\line(-2,-1){20}}

\put(140,50){\line(2,1){20}}

\put(140,50){\line(2,-1){20}}

\put(180,50){\line(-2,1){20}}

\put(180,50){\line(-2,-1){20}}

\put(220,50){\line(2,1){20}}

\put(220,50){\line(2,-1){20}}

\put(260,50){\line(-2,1){20}}

\put(260,50){\line(-2,-1){20}}

\put(65,49){$\ldots$}

\put(285,49){$\ldots$}

\put(195,49){$\ldots$}

\put(40,50){\circle*{3}}

\put(60,50){\circle*{3}}

\put(80,50){\circle*{3}}

\put(100,50){\circle*{3}}

\put(180,50){\circle*{3}}

\put(140,50){\circle*{3}}

\put(220,50){\circle*{3}}

\put(260,50){\circle*{3}}

\put(280,50){\circle*{3}}

\put(300,50){\circle*{3}}

\put(320,50){\circle*{3}}
  \end{picture}
\caption{The structure of $G$}
\label{figure2}

\end{figure}
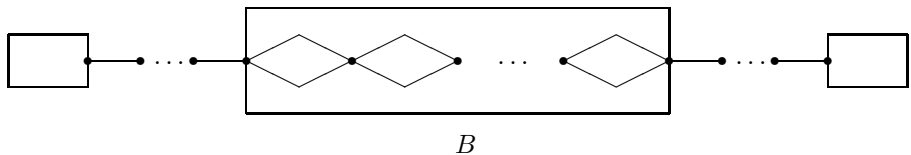

\end{center}

We classify $G$ according to the number of blocks.
If there is no block in $G$, then $G=P_n$. Contradiction!

Suppose that there are at least two blocks in $G$. In this case, take an
essential path which links two blocks. If this essential path has
length $1$, we consider $G'$ obtained by subdividing this
essential edge. If this essential path has length $2$, let $G'=G$.
Let $u-v-w$ be this essential path. Let $G_1$ and $G_2$ be two
induced subgraphs of $G$ so that $G=G_1\cup G_2$ and $G_1\cap
G_2=v$. Note that each block contains at least one non-essential
vertex, which has degree at least $9$.  We have
$$|G_1|\geq 9 \mbox{ and } |G_2|\geq 9.$$
Since $|G_1|+|G_2|=|G'|+1\leq |G|+2$, we have
$$|G_1|< |G| \mbox{ and } |G_2|< |G|.$$
By the minimality of $G$, we have for $i=1,2$,
$$f(G_i)>\sqrt{2}-1.$$
Note $D(G')\leq D(G_1) +D(G_2)$.
\begin{eqnarray*}
  R(G')-R(G_1)-R(G_2) &=& \frac{1}{\sqrt{2d_u}} + \frac{1}{\sqrt{2d_w}}
-\frac{1}{\sqrt{d_u}} - \frac{1}{\sqrt{d_w}}\\
&=& -(1-\frac{1}{\sqrt{2}}) (\frac{1}{\sqrt{d_u}}+\frac{1}{\sqrt{d_w}}) \\
&>&  -(1-\frac{1}{\sqrt{2}}) (\frac{1}{\sqrt{2}}+\frac{1}{\sqrt{2}}) \\
&=& 1-\sqrt{2}.
\end{eqnarray*}
By Lemma \ref{subdivision}, we have
\begin{eqnarray*}
  f(G) &\geq & f(G')\\
       &=&   R(G')-\frac{1}{2}D(G')\\
       &=& f(G_1) + f(G_2)+ R(G')-R(G_1)-R(G_2)\\
       &>&  f(G_1) + f(G_2)+ 1-\sqrt{2}\\
       &>& \sqrt{2}-1 + \sqrt{2}-1 +  1-\sqrt{2}\\
       &=& \sqrt{2}-1.
\end{eqnarray*}
Contradiction!

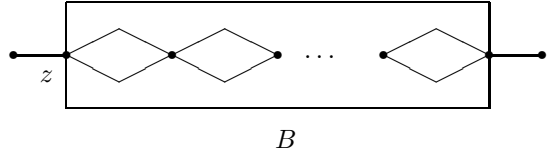
\begin{figure}[hbt]
  \centering
  \begin{picture}(360, 100)






\put(80,50){\line(1,0){20}}

\put(100,30){\line(0,1){40}}

\put(260,30){\line(0,1){40}}

\put(100,30){\line(1,0){160}}

\put(179,15){$B$}

\put(100,70){\line(1,0){160}}

\put(260,50){\line(1,0){20}}






\put(100,50){\line(2,1){20}}

\put(100,50){\line(2,-1){20}}

\put(140,50){\line(-2,1){20}}

\put(140,50){\line(-2,-1){20}}

\put(140,50){\line(2,1){20}}

\put(140,50){\line(2,-1){20}}

\put(180,50){\line(-2,1){20}}

\put(180,50){\line(-2,-1){20}}

\put(220,50){\line(2,1){20}}

\put(220,50){\line(2,-1){20}}

\put(260,50){\line(-2,1){20}}

\put(260,50){\line(-2,-1){20}}






 \put(80,50){\circle*{3}}

\put(100,50){\circle*{3}}

\put(180,50){\circle*{3}}

\put(140,50){\circle*{3}}

\put(220,50){\circle*{3}}

 \put(190,49){$\ldots$}

\put(260,50){\circle*{3}}

\put(280,50){\circle*{3}}


\put(90,40){$z$}

\end{picture}
  \caption{$G$ contains exactly one block with optional
essential edges attached at the end.}
  \label{figure3}
\end{figure}

Now we consider the remaining case: there is exactly one block in $G$
with possible essential edges attached at one or both ends. (See Figure
\ref{figure3}.)

Assume the maximum degree $\Delta$ is achieved at  vertex $v$.
Note that the neighborhood of $v$ can contain at most two essential vertices.
An essential vertex has degree  at least 2 while a non-essiential vertex has degree
at least $9$. Applying Lemma \ref{l0}, we have
\begin{eqnarray*}
  R(G)&\geq& \frac{\sum_{i=1}^n \sqrt{d_i}}{2 \sqrt{\Delta}}\\
&\geq& \frac{\sqrt{d_v}+ \sum_{u\in \Gamma(v)}\sqrt{d_u}}{2 \sqrt{\Delta}}\\
&\geq& \frac{1}{2} + \frac{(\Delta-2) \sqrt{9} + 2 \sqrt{2}}{2 \sqrt{\Delta}}\\
&=& \frac{3}{2}\sqrt{\Delta} -\frac{3-\sqrt{2}}{\sqrt{\Delta}} + \frac{1}{2}.
\end{eqnarray*}

Let $h(x)= \frac{3}{2}\sqrt{x} -\frac{3-\sqrt{2}}{\sqrt{x}} + \frac{1}{2}$.
Note that $h(x)$ is an increasing function on $(0,\infty)$.
Since $G$ contains at least one non-essential vertex, we have $\Delta\geq 9$.

If $D(G)\leq 8$, then we have
$$R(G)\geq h(9) = 4+ \frac{\sqrt{2}}{3}> \frac{D(G)}{2}+\sqrt{2}-1.$$

It remains to show the case $D(G)\geq 9$. In fact, we show the
maximum degree $\Delta$ grows exponentially as a function of $D(G)$.

Pick any path $Q$ (in $G$) of length $D(G)$. Any optional
essential edge(s) is located at the end(s) of $Q$. Let $P$ be the
remaining path after deleting essential edges from $Q$. Let $k$ be
the length of $P$, which is called the length of the block $B$.
Since $D(G)\geq 9$, we have $$k\geq D(G)-2\geq 7.$$

Let $z$ be  an end vertex of $P$.  
For $0\leq i\leq k$, let
$A_i$ be the set of vertices in $B$ of distance $i$ to the
vertex $z$ (see Figure \ref{figure3}).
Let $a_i$ be the minimum degree of
nonessential vertices in $A_i$. If $A_i$ is a single essential
vertex, then define $a_i$ to be infinity.
We have the following two claims.

\noindent
{\bf Claim A:} If $3\leq i\leq k-3$, then we have
\begin{equation}\label{recurrence}
 a_i\geq 2.9(\min\{a_{i-2}, a_{i-1}, a_{i+1}, a_{i+2}\}-1).
\end{equation}

\noindent {\bf Claim B:}  We have $\Delta \geq 1.5+7.4\cdot 2.9^{\lceil (k-6)/4\rceil}$
for $k\geq 7$.

The proofs of these two claims are quite long. We leave these proofs
at the end of this section. Now we use these claims to prove $f(G)>\sqrt{2}-1$. 
For $k\geq 7$, we have
\begin{eqnarray*}
f(G)&=& R(G)-\frac{D(G)}{2}\\
&\geq&  h(\Delta) -\frac{k+2}{2}\\
&\geq&  h(1.5+7.4\cdot 2.9^{\lceil (k-6)/4\rceil}) -\frac{k+2}{2}\\
&>& \sqrt{2}-1.
\end{eqnarray*}
The inequality in last step can be easily verified by Calculus.

The proof of theorem is finished. \hfill $\square$

It remains to prove the two claims.

{\bf Proof of Claim A:} Obviously, (\ref{recurrence}) holds if $a_i$ is
infinity. Suppose there exists  $i$  such that
$$a_i< 2.9(\min\{a_{i-2},
a_{i-1}, a_{i+1}, a_{i+2}\}-1).$$

Let $v$ be the non-essential vertex with degree $a_i$ in $A_i.$
Let $\delta=\min\{a_{i-2}, a_{i-1},$ $ a_{i+1}, a_{i+2}\}.$ The
above inequality implies
\begin{equation}
  \label{eq:5}
d_v<2.9(\delta-1).
\end{equation}

We need show $R(G)>R(G - v)$ to derive the contradiction. It suffices to show that for any
$u\in \Gamma^*(v)$ the weak deletion condition holds.


If $u$ is essential, then $u$ is not connected with any other
essential vertex. We have

\begin{eqnarray*}
&& \hspace*{-1cm} \frac{1}{d_u-1}\sum_{x\in \Gamma(u)\setminus
\{v\}}\frac{1}{\sqrt{d_x^*}}\\
 &\leq& \frac{1}{d_u-1}\sum_{x\in \Gamma(u)\setminus
 \{v\}}\frac{1}{\sqrt{\delta_v-1}} \\
&=& \frac{1}{\sqrt{\delta_v-1}}\\
&<&\frac{2}{\sqrt{d_v}}.
\end{eqnarray*}
At the last step, we applied inequality (\ref{eq:5}).

Otherwise, $u$ can only be adjacent to at most two essential
vertices. Since no two essential vertices are connected, each
non-leaf essential vertex has a degree at least $3$.
Let $y_1$ and $y_2$ be two possible essential vertices.
Noticing that essential vertices are not adjacent, we can orient
the edges of $G\mid_{\Gamma(u)}$ such that directed edges always
leave essential vertices. For  $i\in \{1,2\}$, we have
$$d_{y_i}-\epsilon_{y_i}^u=d_{y_i}\geq 3.$$
For $x\in \Gamma(u)\setminus \{v,y_1,y_2\}$, we apply the bound
$$d_x-\epsilon_{x}^u\geq d_x-1\geq \delta-1.$$
We have
$$\frac{1}{d_u-1}\sum_{x\in \Gamma(u)\setminus
\{v\}}\frac{1}{\sqrt{d_x-\epsilon_x^u}}\leq
\frac{1}{\sqrt{\delta-1}}+\frac{2}{\delta-1}(\frac{1}{\sqrt{3}}-\frac{1}{\sqrt{\delta-1}}).$$

Let
$f(x)=\frac{1}{\sqrt{x}}+\frac{2}{x}(\frac{1}{\sqrt{3}}-\frac{1}{\sqrt{x}})$. Note
that $f(x)$ is decreasing on $(2.5,\infty)$.
Since $d_v\geq 9$, we have $2.5<\frac{d_v}{2.9}<\delta-1$.
Thus, $f(\delta-1)\leq
f(\frac{d_v}{2.9})$. We have
\begin{eqnarray*}
&&\hspace*{-2cm}\frac{1}{d_u-1}\sum_{x\in \Gamma(u)\setminus
\{v\}}\frac{1}{\sqrt{d_x-\epsilon_x^u}}\\
 &\leq& f(\delta-1)\\
&<& f\left(\frac{d_v}{2.9}\right)\\
&<& \frac{2}{\sqrt{d_v}}.
\end{eqnarray*}
The last step can be easily verified by Calculus. \hfill $\square$

{ \bf Proof of Claim B:}
 Let $\{b_i\}$ be the sequence such that $b_i=2.9(b_{i-1}-1)$ and $b_0=9.$
Solving the recurrence equation of the sequence $\{b_i\}$, we have
 \begin{eqnarray*}
 b_i&=& \frac{29}{19}+ \frac{142}{19}\cdot 2.9^i\\
&>& 1.5+7.4\cdot 2.9^i.
\end{eqnarray*}
Since $a_i\geq 9= b_0$ for $0\leq i \leq k$, we get $a_i\geq  b_1$ for $3\leq i\leq k-3$
by applying inequality (\ref{recurrence}).  Applying inequality (\ref{recurrence}) again,
we obtain $a_i\geq  b_2$ for $5\leq i\leq k-5$. Repeatedly apply inequality
(\ref{recurrence}). For each $j$ in $\{1,2,\ldots,\lceil (k-6)/4\rceil  \}$
and each $i$ satisfying $2j+1\leq i \leq k-2j-1$, we have $a_i> b_j$.
Let $j_0=\lceil (k-6)/4\rceil $. Note  $k-4j_0-2\geq 1$. Thus,
both $a_{2j_0+1}$ and $a_{2j_0+2}$ are greater than or equal to $b_{j_0}$.
Note that there is no  essential edge in the block $B$.   We have
$$\Delta\geq \min\{a_{2j_0+1}, a_{2j_0+2}\}\geq b_{j_0}>1.5+7.4\cdot 2.9^{\lceil (k-6/)4\rceil  }.$$
\hfill $\square$

\end{document}